\documentclass{siamart190516}

\usepackage{lipsum}
\usepackage{amsfonts}
\usepackage{graphicx}
\usepackage{amssymb}
\usepackage{algorithm}
\usepackage{algpseudocode}
\ifpdf
  \DeclareGraphicsExtensions{.eps,.pdf,.png,.jpg}
\else
  \DeclareGraphicsExtensions{.eps}
\fi

\newsiamremark{remark}{Remark}
\newsiamremark{hypothesis}{Hypothesis}
\crefname{hypothesis}{Hypothesis}{Hypotheses}
\newsiamthm{claim}{Claim}

\headers{Non-smooth Variable Projection}{T. van Leeuwen and A. Aravkin}

\title{Variable Projection for non-smooth problems\thanks{Submitted to the editors June 29, 2020.}}

\author{Tristan van Leeuwen\thanks{Centrum Wiskunde \& Informatica, Amsterdam, the Netherlands
  (\email{t.van.leeuwen@cwi.nl}).}
\and Aleksandr Aravkin\thanks{Department of Applied Mathematics, University of Washington, Seattle, WA
  (\email{aravkin@uw.edu})}}

\usepackage{amsopn}

\newcommand{\argmin}{\operatornamewithlimits{arg min}}
\newcommand{\prox}{\operatornamewithlimits{prox}}
\newcommand{\ts}[1]{{\textstyle#1}}

\begin{document}
\maketitle
\begin{abstract}
Variable projection solves structured optimization problems
by completely minimizing over a subset of the variables while iterating
over the remaining variables.  Over the last 30 years,
the technique has been widely used, with empirical and theoretical results demonstrating both greater efficacy
and greater stability compared to competing approaches.
Classic examples have exploited closed-form projections and smoothness of the objective function.
We extend the approach to problems that include non-smooth terms, and where
the projection subproblems can only be solved inexactly by iterative methods.
We propose an inexact adaptive algorithm for solving such problems and analyze its computational complexity.
Finally, we show how the theory can be used to design methods for selected problems occurring frequently in machine-learning and inverse problems.
\end{abstract}

% REQUIRED

% REQUIRED

\section{Introduction}
In this paper we consider finite-dimensional \emph{separable} optimization problems of the form
\begin{equation}
\label{bilevel}
\min_{x,y} f(x,y) + r_1(x) + r_2(y),
\end{equation}
where $f$ smoothly couples $(x,y)$ but may be non-convex, while $r_1$ and $r_2$ encode additional constraints or regularizers.
We are particularly interested in the case where $f(x,\cdot) + r_2$ is strongly convex in $y$, so that fast solvers are be available for optimizing over $y$ for fixed $x$. These problems arise any time non-smooth regularization or constraints are used to regularize certain difficult non-linear inverse problems or regression problems. We give three motivating examples below.

\subsection{Motivating examples}

\subsubsection*{Model calibration} Consider the non-linear fitting problem
\begin{equation}
\label{SeparableLeastSquares}
\min_{x,y} \|A(x)y - b\|^2,
\end{equation}
where $A(x)$ defines a linear model with calibration parameters $x$ and $b$ denotes the data. Well-known examples include exponential data-fitting and model-calibration in inverse problems. Adding regularization terms immediately gives a problem of the form \cref{bilevel}.

\subsubsection*{Machine learning}
Trimming is a model-agnostic tool for guarding against outliers. Given any machine learning model that minimizes some objectives $\ell_i$ over a training set
of $m$ datapoints, we introduce $m$ auxiliary parameters $y$ that serve to distinguish inliers from outliers and solve
\[
\min_{x, y} \sum_{i=1}^m y_i \ell_i(x)  \quad \mbox{s.t.} \quad 0 \leq y_i \leq 1 \, \text{for} \, i \in \{1, \ldots, m\},\quad \sum_{i=1}^m y_i = k,
\]
where $k \leq m$ is the number of datapoints we want to fit. The functions $\ell_i$ can capture a wide range of machine learning models.

\subsection*{PDE-constrained optimization}
Many PDE-constrained optimization problems in data-assimilation, inverse problems and optimal control can be cast as
\[
\min_{x,y} \|Py - d\|^2 + \lambda\|A(x)y - q\|^2,
\]
where $y$ denotes the state of the system, $P$ is the sampling operator, $A(x)y = q$ is the discretized PDE with coefficients $y$ and source term $q$, and $\lambda$ is a penalty parameter. Adding regularization terms or changing the data-fidelity term gives a problem of the form \cref{bilevel}.

\subsection{Approach}
The development of specialized algorithms for~\eqref{bilevel} goes back to the classic {\it Variable Projection} (VP) technique for {\em separable non-linear least-squares problems } of the form \eqref{SeparableLeastSquares}, where the matrix-valued map $x \mapsto A(x)$ is smooth and the matrix $A(x)$ has full rank for each $x$. Early work on the topic, notably by \cite{GolubPereyra} has found numerous applications in chemistry, mechanical systems, neural networks, and telecommunications. See the surveys  of~\cite{GolubPereyra2003} and \cite{Osborne2007}, and references therein.

The VP approach is based on eliminating the variable $y$, as for each fixed $x$ we have a closed-form solution
\[
 \overline{y}(x)  = A(x)^\dagger b,
\]
where $A(x)^\dagger$ denotes the Moore-Pensrose pseudo-inverse of $A(x)$. We can thus express \eqref{SeparableLeastSquares} in reduced form as
\begin{equation}\label{eqn:impl}
\min_{x}~ \|(A(x)A(x)^\dagger - I)b\|^2,
\end{equation}
which is a non-linear least-squares problem. Note that $A(x)A(x)^\dagger - I$ is an orthogonal projection onto the null-space of $A(x)^T$; hence the name \emph{variable projection}.

It was shown by \cite{GolubPereyra} that the Jacobian of $A(x)A(x)^\dagger b$ contains only partial derivatives of $A(x)$ w.r.t. $x$ and does not include derivatives of $\overline{y}(x)$ w.r.t. $x$. \cite{ruhe1980algorithms} showed that when the Gauss-Newton method for~\eqref{SeparableLeastSquares} converges superlinearly, so do certain Gauss-Newton variants for~\eqref{eqn:impl}. Numerical practice shows that the latter schemes actually outperform the former on the account of a better conditioning of the reduced problem.

The underlying principle of the VP method is much broader than the class of separable non-linear least squares problems. For example, \cite{Bell2008ADo,AravkinVanLeeuwen2012} consider the class of problems
\begin{equation}
\label{InverseProblemClass}
\min_{x, y}~ f(x,y)\;,
\end{equation}
where $f$ is a $C^2$-smooth function; the classic VP problem~\eqref{SeparableLeastSquares} is a special case of~\eqref{InverseProblemClass}. Although we generally do not have a closed-form expression for $\overline{y}(x)$, we define it as
\[
\overline{y}(x) = \argmin_y f(x,y),
\]
and express \eqref{InverseProblemClass} using the projected function
\begin{equation}
\label{eq:value}
\overline{f}(x) := f(x,\overline{y}(x)).
%\min_x \overline{f}(x) + r_1(x):= f(x,\overline{y}(x)) + r_1(x).
\end{equation}
{\it Projection} in the broader context of~\eqref{InverseProblemClass} refers to {\it epigraphical projection}~\cite{RTRW}, or partial minimization of $y$. Under mild conditions, $\overline{f}(x)$ is $C^2$-smooth as well and its gradient is given by
\begin{equation}
\label{GradientExpression}
\nabla \overline{f}(x) = \Bigl.\nabla_x f(x,y)\Bigr|_{y=\overline{y}(x)},
\end{equation}
i.e., it is the gradient of $f$ w.r.t. $x$, evaluated at $\overline{y}(x)$~\cite{Bell2008ADo}. Again, we do not need to compute any sensitivities of $\overline{y}(x)$ w.r.t. $x$.
This is seen by formally computing the gradient of $\overline{f}$ using the chain-rule:
\begin{equation}
\label{eq:gradfbar}
\nabla\overline{f}(x) = \Bigl(\nabla_x f(x,y) + \nabla_y f(x,y) \cdot \nabla_x\overline{y}(x)\Bigr|_{y=\overline{y}(x)}.
\end{equation}
Since $\overline{y}(x)$ is a minimizer of $f(x,y)$ it satisfies $\nabla_y f(x,\overline{y}(x)) = 0$ and the second term vanishes. Similarly, the Hessian of $\overline{f}$ is the Schur complement of $\nabla_{yy}^2f$ of the full Hessian of $f$ \cite{ruhe1980algorithms}:
\begin{equation}
\label{eq:Hessfbar}
\nabla^2 \overline{f}(x) = \Bigl(\nabla_{xx}^2 f(x,y) - \nabla_{xy}^2 f(x,y) \left(\nabla_{yy}^2 f(x,y)\right)^{-1}\nabla_{yx}^2 f(x,y)\Bigr|_{y=\overline{y}(x)}.
\end{equation}
It follows that a local minimizer, $\overline{x}$, of $\overline{f}$ together with $\overline{y}(\overline{x})$ constitute a local minimizer of $f$. An interlacing property of the eigenvalues of the Schur complement can be used to show that the reduced problem has a smaller condition number than the original problem \cite{Smith1992}.
The expression for the derivative furthermore suggests that we can approximate the gradient of $\overline f$ when $\overline{y}(x)$ is known only approximately by ignoring the second term.

We may extend this approach to solve problems of the form \eqref{bilevel} by including $r_2$ in the computation of $\overline y (x)$ and using an appropriate algorithm to minimize $\overline f + r_1$.
Define the proximity operator for any function $g$ as follows:
\[
\ts{\prox_{\alpha g}}(z) = \argmin_x \ts{\frac{1}{2\alpha}}\|x-z\|^2 + g(x).
\]
where $\alpha>0$ is any scaling factor or stepsize.  We can now view the entire approach as proximal-gradient descent on the projected function $\overline f:$
\begin{equation}
\label{eq:basicExact}
x_{k+1} = \ts{\prox_{\alpha r_1}} \left(x_k - \alpha\nabla \overline f(x_k)\right)
\end{equation}
where $\alpha$ is an appropriate step and $\nabla \overline f$ is computed using \eqref{GradientExpression}. This gives rise to the following protype algorithm \ref{protoalg}.

\begin{algorithm}
\caption{Prototype VP algorithm for solving \eqref{bilevel}}
\label{protoalg}
\begin{algorithmic}
\Require{Initial iterate, $x_0$, Lipschitz constant, $\overline{L}$, of $\nabla\overline{f}$}
\State{$\alpha  = 1/\overline{L}$}
\State{$k = 0$}
\While{not converged}
\State{$\overline{y}_{k+1} = \argmin_{y} f(x_k,y) + r_2(y)$}
\State{$x_{k+1} = \prox_{\alpha r_1} \left(x_k - \alpha \nabla_x f(x_k,\overline{y}_{k+1})\right)$}
\State{$k = k + 1$}
\EndWhile
\end{algorithmic}
\end{algorithm}

Naively, this approach can be applied to non-smooth problems; however, it is not immediately obvious that it is guaranteed to converge. In particular, the resulting reduced objective $\overline{f}$ may not be smooth and hence the gradient formula \eqref{GradientExpression} may not be valid. This is illustrated in the following example.

\subsubsection*{Examples}
\begin{figure}
\centering
\includegraphics[scale=.8]{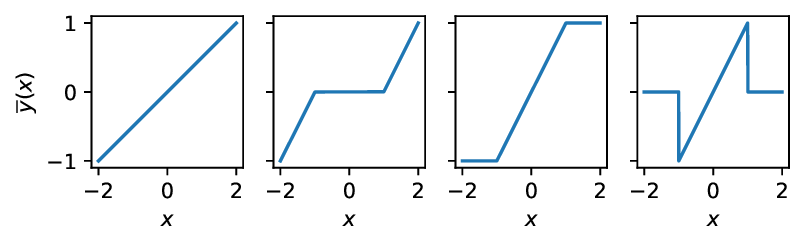}\\
\includegraphics[scale=.8]{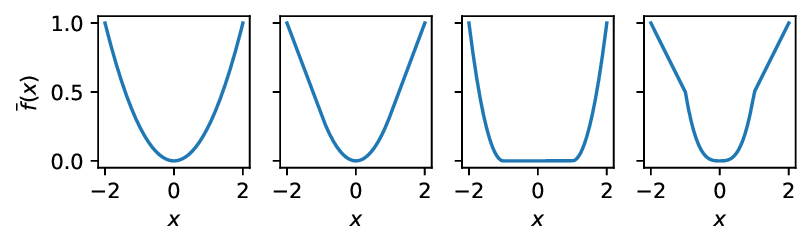}
\label{example1}
\caption{Optimal solution(s) $\overline{y}(x)$ and projected function $\overline{f}(x)$ for functions $F_1$ - $F_4$. For cases 1 - 3, $\overline{y}$ is continuous (but not necessarily smooth), leading to a \emph{smooth} projected function. For $F_4$, $\overline{y}$ is not continuous (at $x = \pm 1$ the solution is not unique),leading to a non-smooth projected function.}
\end{figure}

To illustrate the possibilities and limitation for extending the VP approach to non-smooth functions, consider the following functions
\begin{eqnarray*}
F_1(x,y) &=& \textstyle{\frac{1}{2}}(x - y)^2 + \textstyle{\frac{1}{2}}y^2, \\
F_2(x,y) &=& \textstyle{\frac{1}{2}}(x - y)^2 + |y|,\\
F_3(x,y) &=& \textstyle{\frac{1}{2}}(x - y)^2 + \delta_{[-1,1]}(y), \\
F_4(x,y) &=& |x - y| + x^2 |y|.
\end{eqnarray*}
The corresponding $\overline{y}(x)$ and $\overline{f}(x)$ are shown in figure \cref{example1}. For cases 1 - 3, $\overline{y}$ is continuous (but not smooth), leading to a \emph{smooth} projected function. For $F_4$, $\overline{y}$ is not continuous (at $x = \pm 1$ the solution is not unique), leading to a non-smooth projected function. Continuity of $\overline{y}$ thus appears to be important, rather than smoothness of $F(x,y)$ in $y$. We therefore restrict our attention to problems of the form \cref{bilevel} that are \emph{strongly} convex in $y$ for all $x$ of interest.

The goal of this paper is to extend the variable projection technique to problems of the form~\eqref{bilevel} with \emph{non-smooth} regularization terms, which arise in high-dimensional statistics, signal processing, and many machine learning problems; sparse regularization and simple constraints are frequently used in this setting.

% We note that even simple constraints on $y$ (e.g., $y \in [0,1]^k$) make the problem \eqref{InverseProblemClass} non-smooth.
\subsection{Contributions and outline}
Our contributions are as follows:
\begin{enumerate}
\item  Sufficient conditions under which $\overline{f}$ is smooth and its gradient can be evaluated by \eqref{GradientExpression};
%\item \{Showing how to compute the Lipschitz constant for $\nabla \overline f$, to enable automatic stepsize computation}
\item  Development and analysis of an inexact adaptive version of Algorithm \ref{protoalg} based on inexact evaluations of $\overline{y}(x)$.
\end{enumerate}

In Section~\ref{sec:theory}, we develop derivative formulas for the value function~\eqref{eq:value} and design an inexact version of Algorithm \ref{protoalg}.
In Section~\ref{sec:cases} we present a few case studies. Conclusions complete the paper.

\section{Derivative Formulas and Inexact VP\label{sec:theory}}
\label{section2}
In this section, we present derivative formulas and develop the approaches
briefly described in the introduction. The proofs of the following statements are found in the appendix.
By $\|\cdot\|$ we denote the Euclidean norm. We make the following blanket assumptions on $f$ and $F(x,y) = f(x,y) + r_2(y)$:

\begin{description}
\item[A1] $\nabla_x f$ and $\nabla_y f$ exist and are Lipschitz-continuous for all $(x,y)$:
\begin{eqnarray*}
\|\nabla_x f(x,y) - \nabla_x f(x',y)\| &\leq L_{xx} \|x-x'\|,\\
\|\nabla_x f(x,y) - \nabla_x f(x,y')\| &\leq L_{xy} \|y-y'\|,\\
\|\nabla_y f(x,y) - \nabla_y f(x',y)\| &\leq L_{yx} \|x-x'\|,\\
\|\nabla_y f(x,y) - \nabla_y f(x,y')\| &\leq L_{yy} \|y-y'\|.
\end{eqnarray*}
\item[A2] $F(x,y)$ is $\mu$-strongly convex in $y$ for all $x$.
\end{description}

\subsection{Derivative formulas}

We first establish Lipschitz continuity of the solution $\overline{y}$ w.r.t. $x$  in the following Lemma.
\begin{lemma}
\label{lma:ybar}
Let $F$ satisfy assmptions (A1-A2), then $\overline{y}(x) = \argmin_{y} F(x,y)$ is Lipschitz-continuous as a function of $x$: $$\|\overline{y}(x) - \overline{y}(x')\| \leq (L_{yx}/\mu) \|x - x'\|.$$
\end{lemma}

Next, we establish that the naive derivative-formula \cref{eq:gradfbar} holds under the aforementioned assumptions.
\begin{theorem}
\label{thm:gradient}
Let $F$ satisfy assmptions (A1-A2) and define $\overline{f}(x) = \min_{y} F(x,y)$. The gradient of the projected function is then given by $\nabla \overline{f}(x) = \nabla_x f(x,\overline{y}(x))$, with $\overline{y}(x) = \argmin_y F(x,y)$.
\end{theorem}

Finally, we establish Lipschitz continuity of the gradient of the projected function.

\begin{corollary}
\label{col:lipgrad}
The gradient of $\overline{f}$ is Lipschitz continuous:
$$\|\nabla\overline{f}(x) - \nabla\overline{f}(x')\| \leq \overline{L}\|x - x'\|,$$
where
$$\overline{L} = L_{xx} + L_{xy}L_{yx}/\mu.$$
\end{corollary}

\begin{remark}
Note that the bound in \cref{col:lipgrad} is consistent with the expression for the Hessian in \cref{eq:Hessfbar} in the smooth case.
\end{remark}
% To summarize, these results show that under assumptions (A1-A2) , $\overline{f}$ is Lipschitz smooth with derivative $\nabla \overline{f}(x) = \nabla_x f(x,\overline{y}(x))$.
These results immediately establish the convergence of \cref{protoalg} to a stationary point for a broad class of problems of the form \cref{bilevel}.

\subsection{Inexactness}
\label{sec:inexact}
In many applications, we do not have a closed-form expression for $\overline y$ and it must be computed with an iterative scheme. An obvious choice is to use a proximal gradient method
\[
y_{l+1} = \ts{\prox_{\beta r_2}} \left(y_l - \beta\nabla_y f(x,y_l)\right),
\]
with $\beta \in (0,2/L_{yy})$. The resulting approximation, $\widetilde{y}(x)$, of $\overline y(x)$ yields an approximation of the gradient of $\overline{f}$ with error
\[
\|\nabla_x f(x,\overline{y}(x)) - \nabla_x f(x,\widetilde{y}(x))\| \leq L_{xy} \|\overline{y}(x) - \widetilde{y}(x)\|.
\]
This gives rise to an inexact counterpart of \eqref{eq:basicExact}:
\begin{equation}
\label{eq:basicInexact}
x_{k+1} = \ts{\prox_{\alpha r_1}} \left(x_k - \alpha\nabla_x f(x_k,\widetilde{y}(x_k))\right),
\end{equation}
with $\alpha \in (0,2/\overline{L})$. We denote the error in the gradient as $e_k = \nabla_x f(x_k,\widetilde{y}(x_k)) - \nabla_x f(x_k,\overline{y}(x_k))$. We can immediately bound this error as
\begin{equation}
\label{eq:ebound}
\|e_k\| \leq L_{xy}\|\overline{y}(x_k) - \widetilde{y}(x_k)\|.
\end{equation}
Conditions under which such iterations converge in case $\overline{f}$ is (strongly) convex have been well-studied \cite{devolder2014first,schmidt2011convergence}. The basic gist of these results is that convergence of \cref{eq:basicInexact} can be ensured when the error decays sufficiently fast with $k$.

We did not find any results for general $\overline{f}$ in the literature and therefore establish convergence of the inexact prox-gradient method for general Lipschitz-smooth $\overline{f}$ in the following theorem.

\begin{theorem}[Convergence of inexact proximal gradient - general case]
\label{thm:inexactproxgrad}
Denote
\[
\begin{aligned}
\overline{F}(x) := \overline{f}(x) + r_1(x),
\end{aligned}
\]
and let $\overline{L}$ be the Lipschitz constant of $\nabla\overline{f}$. The iteration
\begin{equation}
x_{k+1} = \text{prox}_{\alpha r_1} \left(x_k - \alpha(\nabla \overline{f}(x_k) + e_k)\right),
\label{eq:basiciter}
\end{equation}
using stepsize $\alpha = 1/\overline{L}$ and with errors obeying $\|e_k\| \leq C \|x_{k+1} - x_k\|$ with $C < \overline{L}/2$ produces iterates for which
\[
\min_{k\in\{0, 1, \ldots, n-1\}} \|x_{k+1} - x_k\| \leq A\sqrt{\frac{\overline{F}(x_0) - \overline{F}(x_*)}{n}},
\]
with $A = \sqrt{\frac{2}{\overline{L}-2C}}$.
\end{theorem}
The proof of this theorem along with supporting Lemma's are presented in Appendix~\ref{sec:inexactProofs}. This result immediately establishes convergence of \cref{eq:basicInexact} for problems satisfying assumptions (A1-A2) when the errors in the gradient obey $\|e_k\| < (\overline{L}/2)\|x_{k+1} - x_k\|$. We discuss how to ensure this in practice in a subsequent section.

Stronger statements on the rate of convergence can be made by making stronger assumptions about $f$.
When $\overline{f}$ is convex (i.e., when $f$ is jointly convex in $(x,y)$), a \emph{sublinear} convergence rate of $\mathcal{O}(1/k)$ can be ensured when $\|e_k\| = \mathcal{O}(1/k^{1+\delta})$ for any $\delta > 0$ {\cite[Prop. 1]{schmidt2011convergence}}. A \emph{linear} convergence for \emph{strongly} convex $\overline{f}$ can be ensured when $\|e_k\| = \mathcal{O}(\gamma^k)$ for any $\gamma < 1$ {\cite[Prop. 3]{schmidt2011convergence}}.

\subsection{Asymptotic complexity}
In this section we analyse the asymptotic complexity of finding an $\overline{\epsilon}$-optimal estimate $\widetilde{x}$ of the solution to \cref{bilevel} which satisfies
\[
|\overline{f}(\widetilde{x}) - \overline{f}(\overline{x})| \leq \overline{\epsilon}.
\]
We measure the complexity in terms of the total number of \emph{inner} iterations required to achieve this.
To usefully analyse the asymptotic complexity we assume that $\overline{f}$ is (strongly) convex.

We first note that we can produce an $\epsilon$-optimal estimate of $\overline{y}$ for which $\|\widetilde{y} - \overline{y}\| \leq \epsilon$ in $\mathcal{O}(\log 1/\epsilon)$ \emph{inner} iterations. This follows directly from linear convergence of the proximal gradient method for strongly convex problems. This also implies that we can compute an approximation of the gradient of $\nabla\overline{f}$ with error bounded by $L_{xy}\epsilon$ in $\mathcal{O}(\log 1/\epsilon)$ inner iterations.

\begin{description}
\item[Convex case] To achieve sublinear convergence of the \emph{outer} iterations for convex $\overline{f}$, we need $\|e_k\| = \mathcal{O}(1/k^{1+\delta})$ for any $\delta > 0$ {\cite[Prop. 1]{schmidt2011convergence}}. To achieve this we need to decrease the inner tolerance at the same rate; $\epsilon_k = \mathcal{O}(1/k^{1+\delta})$. Due to linear convergence of the inner iterations, this requires $\mathcal{O}(\log k)$ inner iterations. A total of $K$ outer iterations thus has an asymptotic complexity of $\mathcal{O}(K \log K)$. Due to sublinear convergence of the outer iterations, we require $K = \mathcal{O}(1/\overline{\epsilon})$ outer iterations, giving a complexity of $\mathcal{O}(1/\overline{\epsilon} \log 1/\overline{\epsilon})$.
\item[Strongly convex case] To achieve linear convergence of the \emph{outer} iterations for strongly convex $\overline{f}$, we need $\|e_k\| = \mathcal{O}(\gamma^k)$ for any $\gamma < 1$ {\cite[Prop. 3]{schmidt2011convergence}}. To achieve this we need to decrease the inner tolerance at the same rate; $\epsilon_k = \mathcal{O}(\gamma^k)$. Due to linear convergence of the inner iterations, this requires $\mathcal{O}(k)$ inner iterations. A total of $K$ outer iterations then require $\mathcal{O}(K^2)$ inner iterations. Due to linear convergence of the outer iterations, we require $K = \mathcal{O}(\log 1/\overline{\epsilon})$ outer iterations, giving an overall complexity of $\mathcal{O}((\log 1/\overline{\epsilon})^2)$.
\end{description}

\subsection{Practical implementation}
A basic proximal gradient method for solving problems of the form \cref{bilevel} is shown in \cref{naivalgo}, where we denote
 \[
 r\left(\begin{bmatrix} x \\ y\end{bmatrix}\right) := r_1(x) + r_2(y).
 \]

\begin{algorithm}
\caption{Prox-gradient for~\eqref{bilevel}}
\label{naivalgo}
\begin{algorithmic}
\Require{Initial iterates, $x_0, y_0$, Lipschitz constant, $L$, of $\nabla f$.}
\State{$\alpha  = 1/L$}
\State{$k=0$}
\While{$\|x_k - x_{k-1}\| + \|y_k - y_{k-1}\|> \overline{\epsilon}$}
\State{$\begin{bmatrix} x_{k+1} \\y_{k+1} \end{bmatrix} = \prox_{\alpha r } \left( \begin{bmatrix}x_k \\ y_k \end{bmatrix} -
\alpha \nabla f\left(\begin{bmatrix} x_k \\y_k \end{bmatrix}\right)\right)$}
%\State{$x_{k+1} = \prox_{\alpha r_1} (x_k - \alpha \nabla_x f(x_k,y_{k}))$}
%\State{$y_{k+1} = \prox_{\beta r_2} (y_k - \alpha \nabla_x f(x_k, y_k)) $}
\State{$k = k + 1$}
\EndWhile
\end{algorithmic}
\end{algorithm}

The conventional VP algorithm with accurate inner solves is given in \cref{exactalgo}.
By \cref{lma:stopping} the stopping criterion for the inner iterations guarantees that at outer iteration $k$, we have $\|y_l - \overline{y}(x_k)\| \leq (2L_{yy}/\mu)\epsilon$ and hence that the error in the gradient is bounded by
\[
\|\nabla f(x_k,\overline{y}(x_k)) - \nabla f(x_k,y_l)\| \leq \left(2L_{xy}L_{yy}/\mu\right) \epsilon.
\]
We use warmstarts for the inner iterations,
which is expected to dramatically reduce the required number of inner iterations, especially when
getting close to the solution.

\begin{algorithm}
\caption{Variable projection for~\eqref{bilevel} with accurate inner solves}
\label{exactalgo}
\begin{algorithmic}
\Require{Initial iterates, $x_0, y_0$, Lipschitz constants $\overline{L},L_{yy}$, tolerances $\epsilon, \overline{\epsilon}$.}
\State{$\alpha  = 1/\overline{L}$}
\State{$\beta  = 1/L_{yy}$}
\State{$k=0$}
\While{$\|x_k - x_{k-1}\| > \overline{\epsilon}$}
\State{$l=0$}
\While{$\|y_l - y_{l-1}\| > \epsilon$}
\State{$y_{l+1} = \prox_{\beta r_2} \left(y_l - \beta \nabla_y f(x_k, y_l)\right) $}
\State{$l = l+1$}
\EndWhile
\State{$y_0 = y_l$}\Comment{Warmstart}
\State{$x_{k+1} = \prox_{\alpha r_1} \left(x_k - \alpha \nabla_x f(x_k,y_{l})\right)$}
\State{$k = k+1$}
\EndWhile
\end{algorithmic}
\end{algorithm}

An inexact version of \cref{exactalgo} can be implemented by specifying a decreasing sequence of tolerances $\{\epsilon_k\}_k$. However, in practice it would be hard to figure out exactly how fast to decrease the tolerance. Inspired by the requirement for convergence of the inexact proximal gradient method we therefore propose a stopping criterion based on the progress in the outer iterations. \cref{thm:inexactproxgrad} requires that the error in the gradient $\|e_k\|$ is bounded by $(\overline{L}/2)\|x_{k+1} - x_k\|$. We can guarantee this by combining various bounds. By \cref{eq:ebound} and \cref{lma:stopping} we get
\[
\|y_l - y_{l-1}\| \leq \frac{\overline{L}\mu}{2L_{yy}L_{xy}}\|x_{k+1} - x_k\|,
\]
which guarantees the required bound on the error. Note that $x_{k+1}$ implicitly depends on $y_l$ and $x_k$ through $x_{k+1} = \prox_{\alpha r_1} \left(x_k - \alpha \nabla_x f(x_k,y_{l})\right)$.

A practical implementation of the inexact method is shown in \cref{adaptalgo}. In practice the constants $\overline{L}, L_{yy}, L_{xy}$ and $\mu$ may difficult to estimate. Moreover, the bounds may be very loose. We therefore introduce a parameter, $\rho$, to use in the stopping criterion instead.

In particular settings, the structure of the inner problem in $y$ can be exploited to achieve superlinear convergence of the inner iterations. A further improvement could be to use a line search or an accelerated proximal gradient method for the outer iterations. However, accelerated proximal gradient methods are generally more sensitive to errors or require more information on the function, as pointed out by \cite{schmidt2011convergence}.

\begin{algorithm}
\caption{Algorithm for~\eqref{bilevel} with adaptive tolerance for the inner solves}
\label{adaptalgo}
\begin{algorithmic}
\Require{Initial iterates, $x_0, y_0$, Lipschitz constants $\overline{L},L_{yy}$, parameter $\rho$}
\State{$\alpha  = 1/\overline{L}$}
\State{$\beta  = 1/L_{yy}$}
\State{$k=0$}
\While{$\|x_k - x_{k-1}\| > \overline{\epsilon}$}
\State{$x_{k+1} = x_k$}
\State{$l=0$}
\While{$\|y_l - y_{l-1}\| > \rho \|x_{k+1} - x_k\|$}
\State{$x_{k+1} = \prox_{\alpha r_1} \left(x_k - \alpha \nabla_x f(x_k,y_{l})\right)$}\Comment{Prospective update of $x_k$}
\State{$y_{l+1} = \prox_{\beta r_2} \left(y_l - \beta \nabla_y f(x_k, y_l)\right) $}\Comment{Update of $y_l$}
\State{$l = l+1$}
\EndWhile
\State{$y_0 = y_l$}\Comment{Warm start}
\State{$k = k+1$}
\EndWhile
\end{algorithmic}
\end{algorithm}

\clearpage
\section{Case studies}
\label{sec:cases}

\subsection{Reproducibility}
The Python code used to conduct the numerical experiments is available at \url{https://github.com/TristanvanLeeuwen/VarProNS}. We implemented \cref{naivalgo} (hereafter referred to as the \emph{joint} approach), \cref{exactalgo} (\emph{VP}) and \cref{adaptalgo} (\emph{adaptive VP}). In stead of the absolute tolerance for the inner iterations in \cref{exactalgo} we use a relative tolerance and stop the inner iterations in \cref{exactalgo} when $\|y_{l} - y_{l-1}\| \leq \epsilon \|y_{l-1}\|$. For the outer iterations we use a stopping criterion on the function value. The Lipschitz constants $L,\overline{L}$, $L_{yy}$ and other algorithmic parameters are specified for each example. We show results for several values of $\rho$ to investigate the sensitivity of the results to this parameter.

Aside from the example-specific results, we report the convergence history (in terms of the value of the objective) as a function of the number of outer iterations and the total number of iterations. The latter is used as an indication of the computational cost, and will be referred to as \emph{cost}.

\subsection{Exponential data-fitting}
We begin with the general class of \emph{exponential data-fitting} problems -- one of the prime applications of variable projection \cite{Pereyra2012}.
The general formulation of these problems assumes a model of the form
\[
{d}_i = \sum_{j=1}^{n} y_j\exp(-\phi_{ij}(x)),\quad i = 1,2, \ldots, m,
\]
where $y\in\mathbb{R}^n$ are unknown weights, $d\in\mathbb{R}^m$ are the measurements, and $\phi_{ij}$ are given functions that depend on an unknown parameter $x\in\mathbb{R}^n$. Some examples of this class are given in table \ref{table:expexample}.

\begin{table}
\caption{Some examples of exponential data-fitting in  applications.}
\label{table:expexample}
\begin{tabular}{p{6cm}|p{2cm}}
application                            & $\phi_{ij}$                   \\
\hline
pharmaco-kinetic modelling             & $x_i t_j$                     \\
multiple signal classification         & $\imath x_i h_j$              \\
radial basis function interpolation    & $\alpha_i^2\|x_i - \xi_j\|^2$ \\
\end{tabular}
\end{table}

The exponential data-fitting problem is traditionally formulated as a least-squares problem
\[
\min_{x,y} \|A(x)y - d\|^2.
\]
In many applications, however, it is natural to include regularization terms and/or use another data fidelity term. For example, \cite{Cornelio2012,Shearer2013} consider positivity constraints $y_i \geq 0$. Another common regularization that enforces sparsity of $y$ is $r_2(\cdot) = \lambda \|\cdot\|_1$. This is useful in cases where the system is over-parametrized and we are looking for a fit of the data with as few components as possible. Even with regularization, we expect a highly ill-posed problem where many parameter combinations lead to a nearly identical data-fit.

For the conditions of the previous theorems to hold we require $A(x)^T\!A(x)$ to be invertible for all $x$. Note that this would require $x_i \not= x_j$ for $i \not= j$. As long as we initialize the $x_i$ to be different, we don't expect any problems. Alternatively, a small quadratic term in $y$ could be added to ensure that $f(x,y) = \|A(x)y - d\|^2 + \delta \|y\|^2$ is strongly convex in $y$.

\subsubsection{Example 1}
\begin{figure}
\centering
\includegraphics[scale=.5]{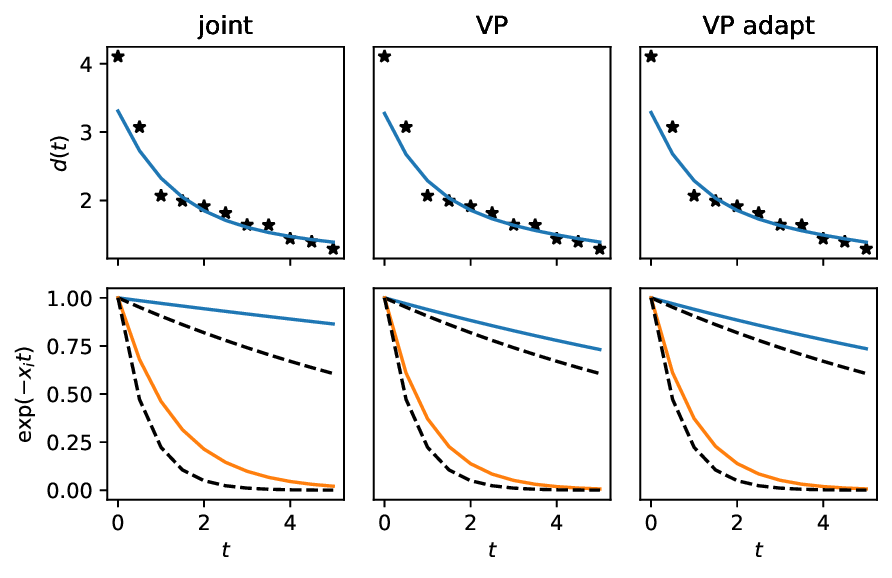}
\label{fig:expfit1b}
\caption{Top: data and fitted data resulting from the various methods. Bottom: individual components $\exp(-x_i t)$ resulting from the optimization, as well as the ground-truth components (black dashed line)}
\end{figure}

\begin{figure}
\centering
\includegraphics[scale=.5]{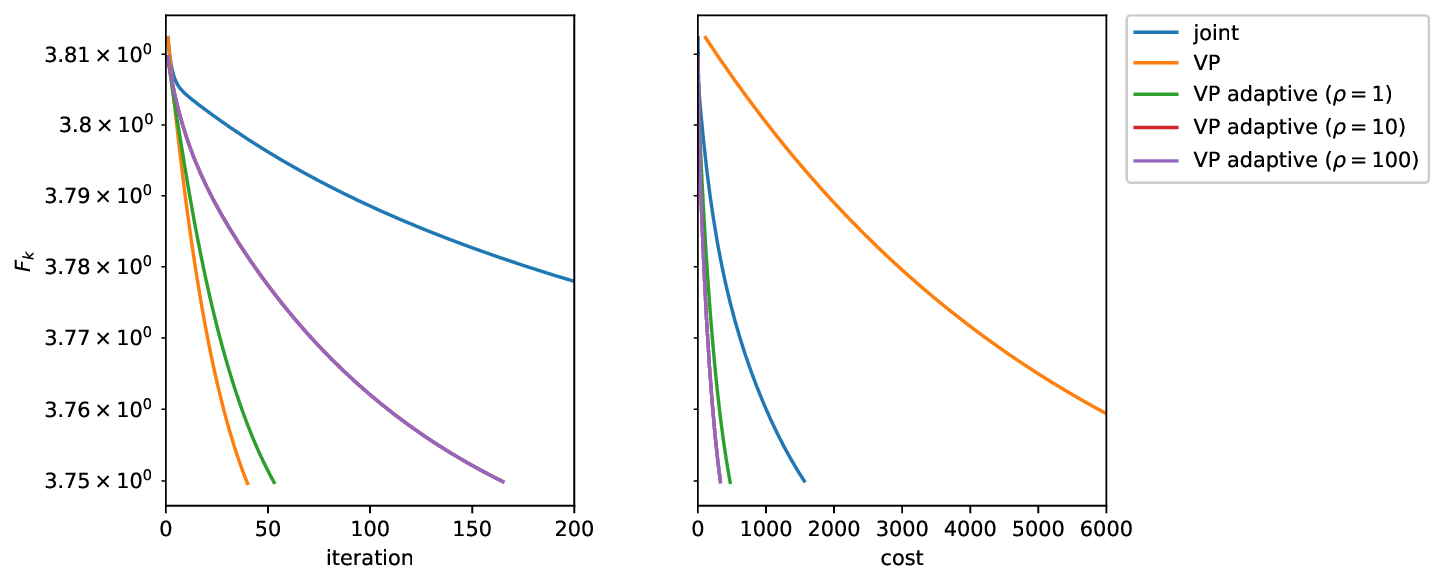}
\label{fig:expfit1a}
\caption{Convergence history and computational cost for \cref{naivalgo}, \cref{protoalg}, \cref{adaptalgo}. Note that the inexact VP algorithm retains the same favourable rate of convergence of the full VP algorithm while being significantly cheaper than the two alternatives.}
\end{figure}

We model the data via a normal distribution with mean $\overline{d}_i = \sum_{j=1}^n y_j \exp(-t_i x_j)$ and variance $\sigma$.
To generate the data, we set $m = 11$, $n = 2$, $t_i = (i-1)/2$, $x = (0.1, 1.5)$, $y = (2.1,1.9)$ and $\sigma=0.1$. Since the number of components, $n$, is not typically known in practice, we attempt the fit the data using $n=5$ components and use an $\ell_1$-regularization term with $\lambda = 1$ (chosen by trial and error) to find a parsimonious solution. The variational problem is then given by
\[
\min_{x,y} \|A(x)y - d\|^2 + \lambda \|y\|_1 \quad \text{s.t.} \quad y_i \geq 0, \,\text{for}\, i = 1,2, \ldots, n.
\]
Here, $A(x)$ is an $m \times n$ matrix with elements $\exp(-t_i x_j)$. The Lipschitz-constants of $f$ and $\overline{f}$ are set to $L = \overline{L} = 1\cdot 10^3$ (selected by trial and error). For the inner solves we let $L_{yy} = \|A(x)\|^2$ (for the true $x$), $\epsilon = 10^{-6}$ for VP and $\rho \in \{1, 10, 100\}$ for adaptive VP. With these settings we run \cref{naivalgo} (joint), \cref{exactalgo} (VP) and \cref{adaptalgo} (adaptive VP). The results are shown in \cref{fig:expfit1b} and \cref{fig:expfit1a}. The VP approach gives a slightly better recovery of the modes $\exp(-x_i t)$ while achieving the same level of data-fit as the joint approach. This indicates that the ill-posedness of the problem is effectively dealt with by the VP approach. The VP approach converges much faster than the joint approach. While the exact VP approach is much more expensive than the joint approach, the adaptive VP method is much cheaper than the joint approach.

\subsubsection{Example 2}
\begin{figure}
\centering
\includegraphics[scale=.5]{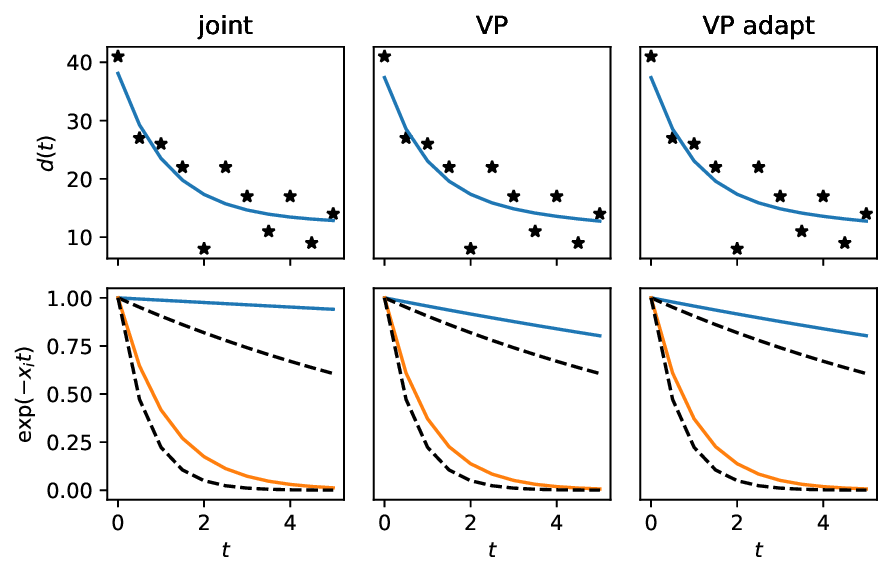}
\label{fig:expfit2b}
\caption{Top: data and fitted data resulting from the various methods. Bottom: individual components $\exp(-x_i t)$ resulting from the optimization, as well as the ground-truth components (black dashed line)}
\end{figure}

\begin{figure}
\centering
\includegraphics[scale=.5]{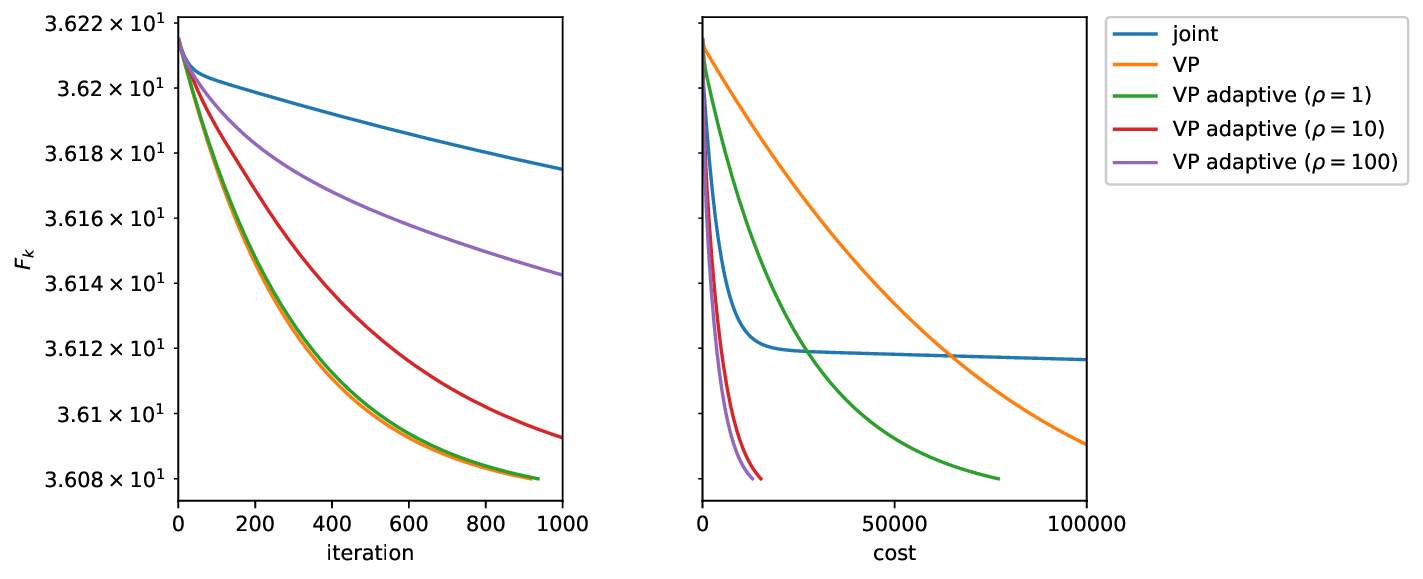}
\label{fig:expfit2a}
\caption{Convergence history and computational cost for \cref{naivalgo}, \cref{exactalgo} and \cref{adaptalgo}. Note that the inexact VP algorithm retains the same favorable rate of convergence of the full VP algorithm while being significantly cheaper than the two alternatives.}
\end{figure}

Here, we model the data via a Poisson distribution with parameter $\overline{d}_i = \sum_{j=1}^n y_j \exp(-t_i x_j)$. We set $m = 11$, $n = 2$, $t_i = (i-1)/2$, $x = (0.1, 1.5)$, $y = (21,19)$. The other settings are exactly the same as the previous example.
The variational problem is now given by
\[
\min_{x,y} \ell(A(x)y,d) + \lambda \|y\|_1 \quad \text{s.t.} \quad y_i \geq 0, \,\text{for}\, i = 1, \ldots, n,
\]
where $\ell$ denotes the Poisson log-likelihood function. The Lipschitz-constants of $f$ and $\overline{f}$ are set to $L = \overline{L} = 5\cdot 10^4$ (selected by trial and error). For the inner solves we let $L_{yy} = \|A(x)\|^2$ (for the true $x$), $\epsilon = 10^{-6}$ for inexact VP and $\rho \in \{1,10,100\}$ for adaptive VP. With these settings we run \cref{naivalgo} (joint), \cref{exactalgo} (VP) and \cref{adaptalgo} (adaptive VP). The results are shown in \cref{fig:expfit2a} and \cref{fig:expfit2b}. Here again the VP approach gives a slightly better recovery of the modes $\exp(-x_i t)$ while achieving the same level of data-fit as the joint approach. This indicated the ill-posedness of the problem that is effectively dealt with by the VP approach. The VP approach converges much faster than the joint approach. The adaptive VP method is (much) cheaper than the joint approach and even the regular VP approach eventually beats the joint approach.

%%%%%%%%%%%%%%%%%%%%%%%%%%%%%%%%%%%%%%%%%%%%%%%%%%%%%%%%%%%%%%%%%%%%%%%%%%%%%%%%%%%%%%%%%

\subsection{Trimmed Robust Formulations in Machine Learning}\label{subsec:trimmed}

Many formulations in high-dimensional regression, machine learning,
and statistical inference can be formulated as minimization problems
\[
\min_x \sum_{i=1}^n \ell_i(x),
\]
where the training set comprises $n$ examples and $\ell_i$ is the error or negative log-likelihood corresponding to the $i$th training point. This approach can be made robust to perturbations of input data
(for example, incorrect features, gross outliers, or flipped labels) using a {\it trimming approach}.
The idea, first proposed by \cite{rousseeuw1984least} in the context of least squares fitting,
is to minimize the $k \leq n$ best residuals. The general trimmed approach, formulated and studied by \cite{yang2018general}, considers the equivalent formulation
\[
\min_{x,y} \sum_{i=1}^n y_i \ell_i(x), \quad y \in \hat \Delta_k,
\]
where $\hat \Delta_k := \{y \in [0,1]^n \, : \, \quad 1^Ty = k\}$ denotes the {\it capped simplex} and admits an efficient projection \cite{wang2015projection,aravkin2016smart}.
Jointly solving for $(x,y)$ selects the $k$ inliers as the model $x$ is fit.
Indeed, the reader can check that the solution in $y$ for fixed $x$ selects the smallest $k$ terms $f_i$. The problem therefore looks like a good candidate for VP, but the projected function $\overline f(x)$ is nonsmooth, because the solution for $y$ may not be unique\footnote{Consider, for example, a case where a number of residuals have the same value.}. However, the smoothed formulation
\begin{equation}
\label{eq:smoothTrim}
\min_{x,y} \sum_{i=1}^ny_i \ell_i(x) + \textstyle{\frac{\delta}{2}}\|y\|^2,\quad y \in \hat \Delta_k,
\end{equation}
does lead to a differentiable $\overline f_{\delta}(x)$ as the problem is now strongly convex in $y$.

To illustrate the trimming method and the potential benefit of the VP approach we consider the following stylized example.

\subsubsection{Example 1}
\begin{figure}
\centering
\includegraphics[scale=.5]{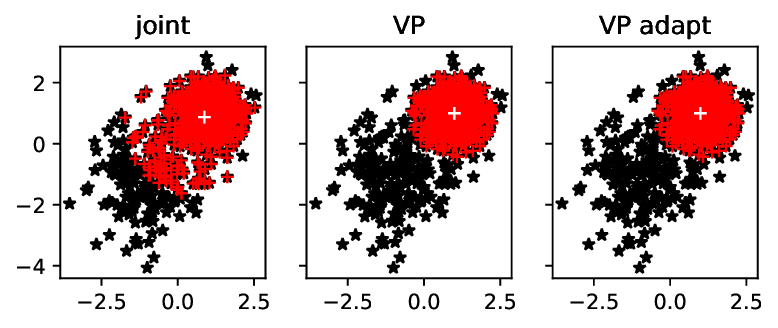}
\label{fig:trimming1b}
\caption{Samples $d_i$ (black), selected inliers (red) and estimated mean (white) for each method.}
\end{figure}

\begin{figure}
\centering
\includegraphics[scale=.4]{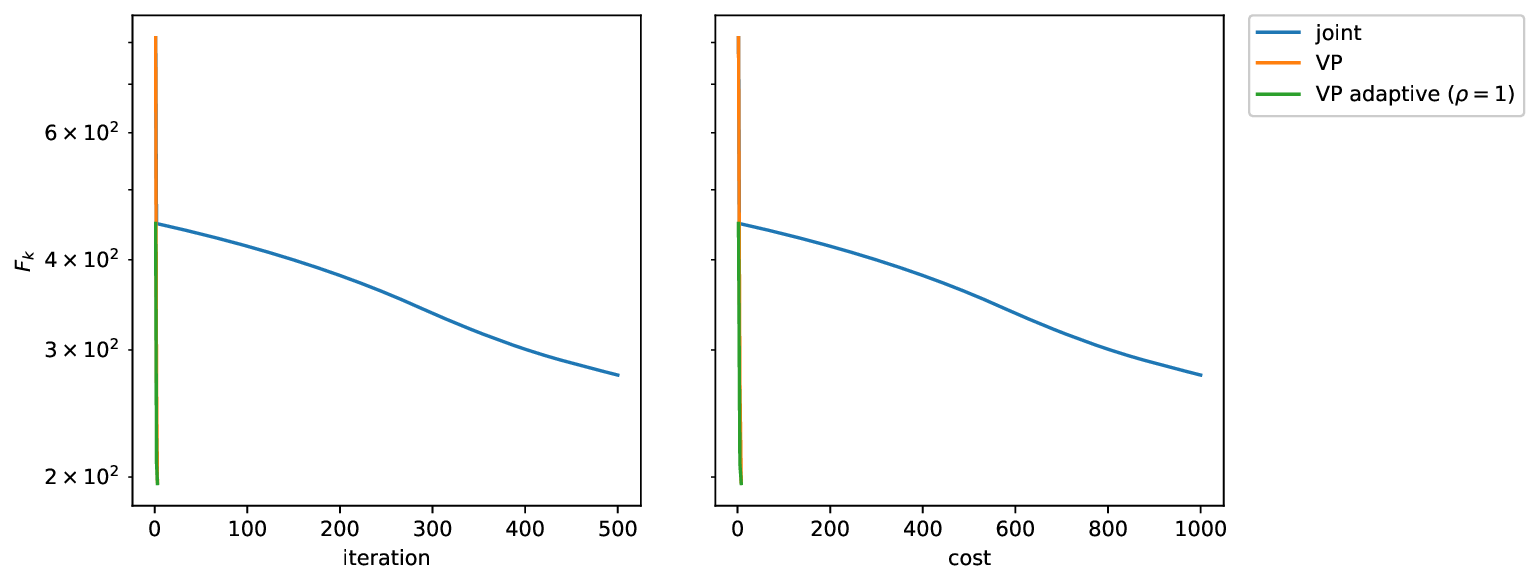}
\includegraphics[scale=.4]{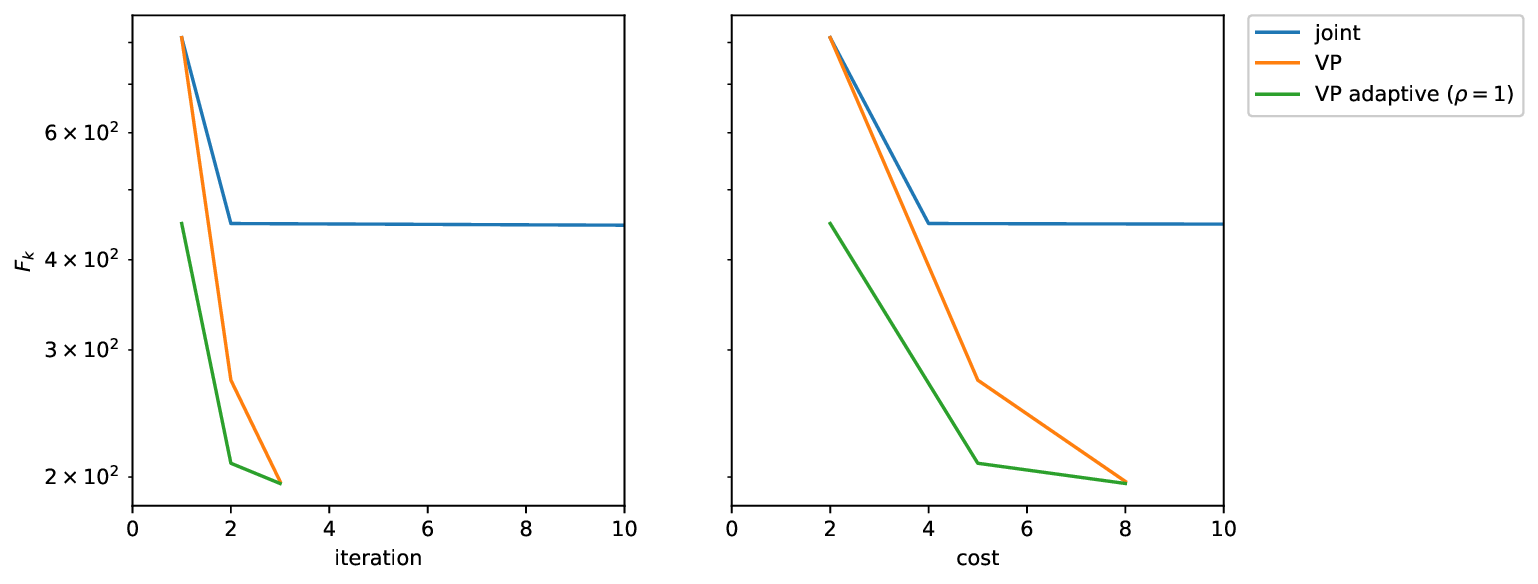}
\label{fig:trimming1a}
\caption{Convergence history and computatonal cost for\cref{naivalgo}, \cref{protoalg}, \cref{adaptalgo}. The adaptive VP algorithm retains the same favorable rate of convergence of the full VP algorithm while being significantly cheaper than the two alternatives. The bottom plots show more clearly what happens in the first few iterations.}
\end{figure}

We aim to compute the (trimmed) mean of a set of samples $\{d_i\}_{i=1}^n$ by setting $\ell_i(x) = \frac{1}{2}\|x - d_i\|^2$. We generate data by sampling generating a total of $n = 1000$ samples $d_i \in \mathbb{R}^{2}$; $k = 200$ are drawn from a normal distribution with mean $(-1,-1)$ and variance $1$ while the remaining $800$ are drawn from a normal distribution with mean $(1,1)$ and variance $0.5$. We let $\delta = 1\cdot 10^{-3}$, $L = \overline{L} = k$ and $L_{yy} = \delta$, $\epsilon = 10^{-6}$ for VP and $\rho = 1$ for adaptive VP. The results for \cref{naivalgo} (joint), \cref{protoalg} (VP), \cref{adaptalgo} (adaptive VP) are shown in figures \cref{fig:trimming1b} and \cref{fig:trimming1a}. We note that the VP approach gives a better recovery of the inliers. Moreover, the VP approach converges much faster than the joint approach. As the inner problem is a simple quadratic, both the exact and adaptive VP approach have essentially the same computational cost.

\subsection{Tomography}\label{sec:tomo}

In computed tomography (CT) the data consists of a set of line integrals of an unknown function
\[
d_{ij} = \int_{\ell_{ij}} u(t) \mathrm{d}t,
\]
where $\ell_{ij}$ represents a line with angle $\theta_i$ and offset $o_j$. Collecting measurements for angles $\{\theta_i\}_{i=1}^k$, offsets $\{o_j\}_{j=1}^m$ and representing the function $u$ as a piecewise continuous function on a grid of $n\times n$ rectangular pixels leads to a system of linear equations $d = Ay$, with $A \in \mathbb{R}^{m\cdot k \times n^2}$, with $y_i$ representing the gray value of the image in the $i^{\text{th}}$ pixel. In many applications, the angles and absolute offsets are not known exactly due to calibration issues. To model this we include additional parameters $x_i = (\Delta \theta_i, \Delta s_i)$ for each angle. This leads to the following signal model
\[
d = A(x)y + e,
\]
where $x \in \mathbb{R}^{2k}$ contains the calibration parameters and $e$ represents measurement noise\footnote{In low-dose applications, a Poisson noise model is more realistic.}. We can ensure that $A(x)$ depends smoothly on $x$ by using higher-order interpolation \cite{van2017automatic}. Many alternative methods have been proposed to solve the calibration problem \cite{van2017automatic,Austin2019} and we do not claim that the proposed VP approach is superior; we merely use this problem as an example to compare the joint, VP and adaptive VP approaches amongst each other.

\subsubsection{Example 1}
We let $n = 50$, and take $m = 50$ offsets regularly sampled in $[-1.5,1.5]$ and $k=50$ angles regularly sampled in $[0,2\pi]$. This does not guarantee an invertible Hessian $A(x)^T\!A(x)$ (e.g., strongly convex $f(x,\cdot)$) for all $x$, however, we will see that this does not lead to numerical issues in this particular example. As a safeguard, one could add a small quadratic regularizer $\delta \|y\|^2$ to the objective.

The data are generated with an additional random perturbation on the offset and angles (both are normally distributed with zero mean and variance $0.5$) and measurement noise (normally distributed with zero mean and variance $1$).

\begin{figure}
\centering
\includegraphics[scale=.8]{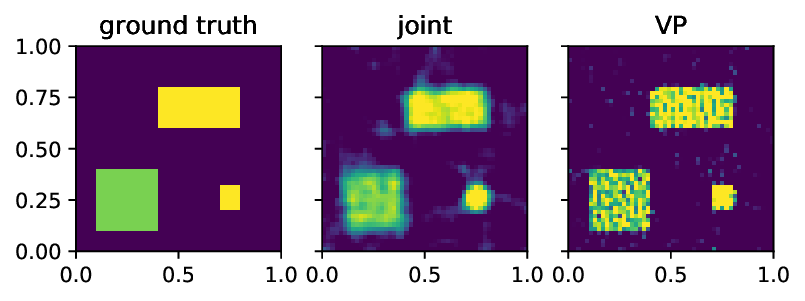}\\
\includegraphics[scale=.4]{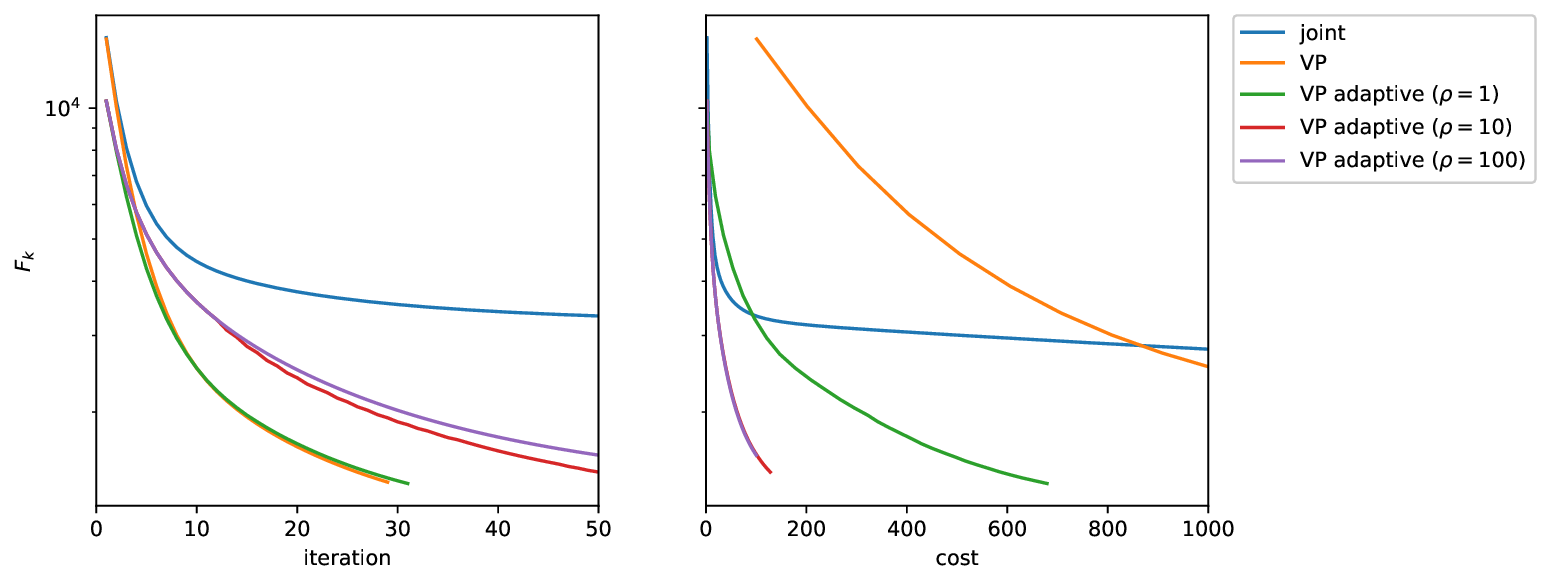}\\
\label{fig:tomo1}
\caption{Top: ground truth and reconstructed images. Bottom: convergence history and computational cost for \cref{naivalgo}, \cref{exactalgo} and \cref{adaptalgo}. The inexact VP algorithm retains the same favorable rate of convergence of the full VP algorithm while being significantly cheaper than the two alternatives.}
\end{figure}

To regularize the problem we include positivity constraints:
\[
\min_{x,y} \|A(x)y - d\|^2 \quad \text{s.t.} \quad y_i \geq 0.
\]
For the optimization we let $L = \overline{L} = 1\cdot 10^6$ and $L_{yy} = \|A(x)\|^2$ (for the true $x$), $\epsilon = 10^{-6}$ for VP and $\rho \in \{1,10,100\}$ for adaptive VP. With these settings we run \cref{naivalgo} (joint), \cref{exactalgo} (VP) and \cref{adaptalgo} (adaptive VP). The results are shown in \cref{fig:tomo1}. The VP approach gives a much better reconstruction than the joint approach, owing to the much faster convergence. The adaptive VP methods are cheaper than the joint method, again because they converge much faster. The adaptive VP approach, in turn, is much cheaper than the exact VP approach.

\subsubsection{Example 2}
\begin{figure}
\centering
\includegraphics[scale=.8]{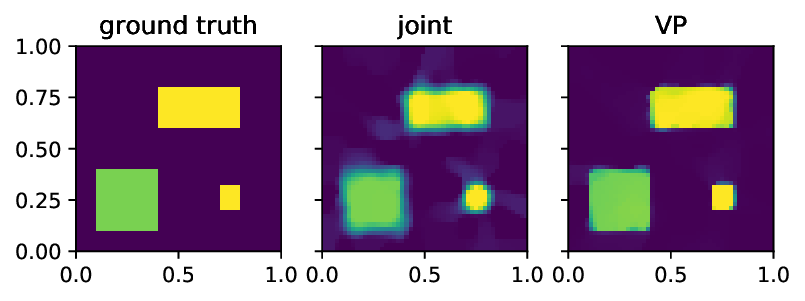}\\
\includegraphics[scale=.4]{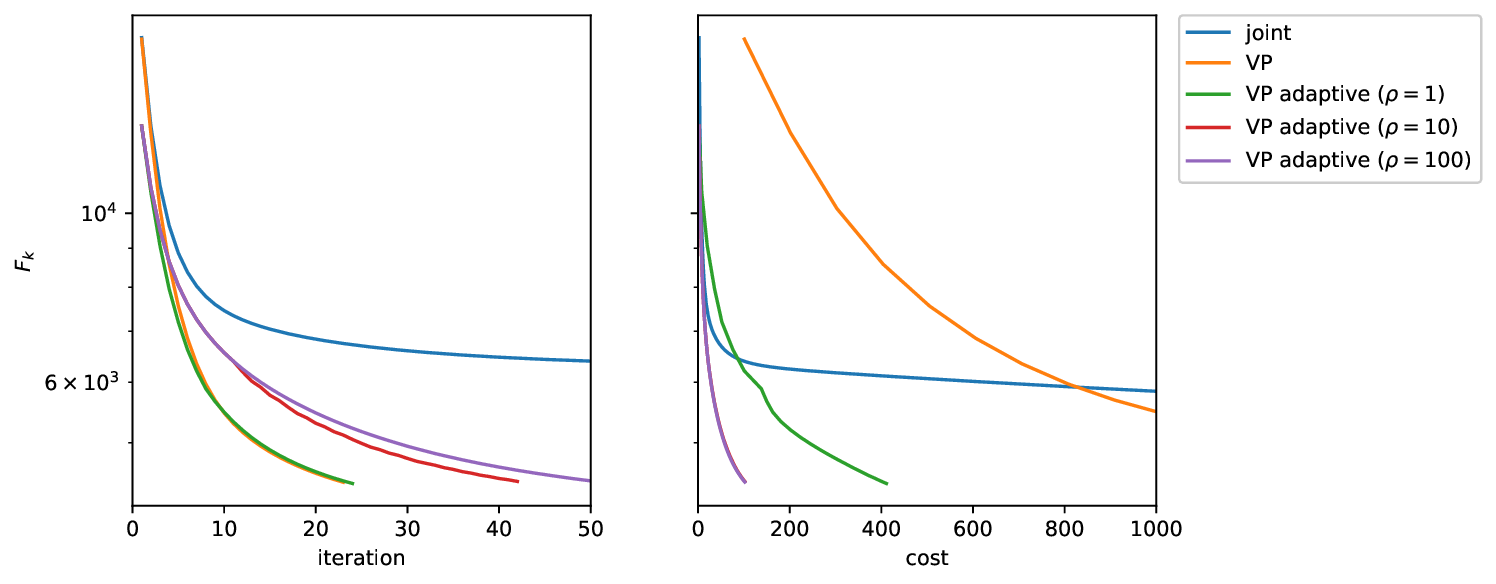}\\
\label{fig:tomo2}
\caption{Top: ground truth and reconstructed images. Bottom: convergence history and computational cost for \cref{naivalgo}, \cref{exactalgo} and \cref{adaptalgo}. Note that the inexact VP algorithm retains the same favorable rate of convergence of the full VP algorithm while being significantly cheaper than the two alternatives.}
\end{figure}

The settings are the same as in the previous example, except that we use a TV-regularization term:
\[
\min_{x,y} \|A(x)y - d\|^2 + \lambda \text{TV}(y),
\]
with
\[
\text{TV}(y) = \sum_{i} \sqrt{(D_1 y)_i^2 + (D_2 y)_i^2},
\]
where $D_k y$ returns a finite-difference approximation of the first derivative of the image $y$ in the $k^{\text{th}}$ direction. We let $\lambda = 2\cdot 10^1$ (chosen by trial and error). The results are shown in \cref{fig:tomo2}. The VP approach gives a much better reconstruction than the joint approach, owing to the much faster convergence. Both exact and adaptive VP are cheaper than the joint method, again because they converge much faster. The adaptive VP approach, in turn, is much cheaper than the exact VP approach.

\clearpage
\section{Conclusions\label{sec:conclusions}}
Variable projection has been successfully used in a variety of
contexts; the popularity of the approach is largely due to its
superior numerical performance when compared to joint optimization schemes.
In this paper, we extent its use to wide class of non-smooth and constrained problems occurring in various applications. In particular, we give sufficient conditions for the applicability of an inner-outer proximal gradient method. We also propose an inexact algorithm with an adaptive stopping criterion for the inner iterations and show that it retains the same convergence rate as the exact algorithm.
Numerical examples on a wide range of nonsmooth applications show that: \emph{i)} the variable projection approach leads to faster convergence than the joint approach, and \emph{ii)} the adaptive variable projection method outperforms both the joint method and exact variable projection method in terms of computational cost. The adaptive method includes a parameter, $\rho$, that controls the stopping tolerance for the inner iterations. In the numerical experiments we observe that a larger value for $\rho$ leads to a smaller number of outer iterations. There appears to be little danger of setting $\rho$ too large; even when varying it over two orders of magnitude the adaptive VP method is consistently faster than the exact VP method. As $\rho \rightarrow 0$, the adaptive method coincides with the exact method. Heuristics could be developed to set this parameter automatically, but since this will be highly application-specific, it is outside the scope of this paper.

\clearpage
\appendix

\section{Proof of \cref{lma:ybar}}
\begin{proof}
Denote $F(x,y) = f(x,y) + r_2(y)$. For ease of notation we fix $x,x'$ and denote the corresponding (unique) optimal solutions by $\overline{y}$ and $\overline{y}'$ respectively. These optimal solutions are implicitly defined by the first order optimality conditions $0 \in \partial_y F(x,\overline{y})$ and $0 \in \partial_y F(x',\overline{y}')$, or equivalently
\begin{equation}
\label{eq:opt}
-\nabla_y f(x,\overline{y}) \in \partial r_2(\overline{y}), \quad
-\nabla_y f(x',\overline{y}') \in \partial r_2(\overline{y}').
\end{equation}
We start from strong convexity of $F$. For \emph{any} $\widetilde{x}$ we have
\begin{equation}
\label{eq:strongcvx}
\mu \|\overline{y}' - \overline{y}\|^2 \leq \langle g'_{\widetilde{x}} - g_{\widetilde{x}}, \overline{y}' - \overline{y}\rangle,
\end{equation}
with $g_{\widetilde{x}} \in \partial_y F(\widetilde{x},\overline{y})$ and $g'_{\widetilde{x}} \in \partial_{y} F(\widetilde{x}, \overline{y}')$. We note that $\partial_y F(\widetilde{x},\overline{y}) = \{\nabla_y f(\widetilde{x},\overline{y}) + h \, |\, h \in \partial r_2(\overline{y})\}$ so we can write
\[
\begin{aligned}
g_{\widetilde{x}} = \nabla_y f(\widetilde{x},\overline{y}) + h, \quad h \in \partial r_2(\overline{y}), \\
g_{\widetilde{x}}' = \nabla_y f(\widetilde{x},\overline{y}') + h', \quad h' \in \partial r_2(\overline{y}'),
\end{aligned}
\]
with~\eqref{eq:strongcvx} holding for any choice of $h$ and $h'$.
Using~\eqref{eq:opt}, we make the particular choices  $h = -\nabla_y f(x, \overline{y}) \in \partial r_2(\overline{y})$ and $h' = -\nabla_y f(x', \overline{y}') \in \partial r_2(\overline{y}')$. From~\eqref{eq:strongcvx} we now have
\[
\mu \|\overline{y}' - \overline{y}\|^2 \leq \langle \nabla_y f(\widetilde{x},\overline{y}') - \nabla_y f(x',\overline{y}'), \overline{y}' - \overline{y}\rangle + \langle \nabla_y f(x,\overline{y}) - \nabla_y f(\widetilde{x},\overline{y}), \overline{y}' - \overline{y}\rangle.
\]
Setting $\widetilde{x} = x$, we get
\[
\mu \|\overline{y}' - \overline{y}\|^2 \leq \langle \nabla_y f(x,\overline{y}') - \nabla_y f(x',\overline{y}'), \overline{y}' - \overline{y}\rangle.
\]
Finally, using Cauchy-Schwartz and Lipschitz-smoothness of $f$ in $y$ we have
\[
\|\overline{y}' - \overline{y}\| \leq (L_{xy}/\mu) \|x' - x\|.
\]
\end{proof}

\section{Proof of \cref{thm:gradient}}
\begin{proof}
We set out to show that
\begin{equation}
\label{eq:limres}
\lim_{\|e\|\rightarrow 0} \frac{|\overline{f}(x + e) - \overline{f}(x) - \nabla_x f(x,\overline{y}(x))\cdot e|}{\|e\|} = 0,
\end{equation}
which would confirm that $\nabla_x f(x,\overline{y}(x))$ is indeed the gradient of $\overline{f}(x)$.

Using the definition of $\overline{f}$ we rewrite~\eqref{eq:limres} as
\[
\lim_{\|e\|\rightarrow 0} \frac{|F(x + e,\overline{y}(x + e)) - F(x,\overline{y}(x)) - \nabla_x f(x,\overline{y}(x))\cdot e|}{\|e\|}.
\]
Writing $F(x + e, \overline{y}(x + e)) = F(x, \overline{y}(x + e)) + \nabla_x f(\xi,\overline{y}(x + e))$ with $\xi = x + t e$ for $t\in [0,1]$, we get
\begin{equation}
\label{eq:limit}
\lim_{\|e\|\rightarrow 0} \frac{|F(x,\overline{y}(x + e)) +
\nabla_x f(\xi,\overline{y}(x + e))\cdot e - F(x,\overline{y}(x)) - \nabla_x f(x,\overline{y}(x))\cdot e|}{\|e\|}.
\end{equation}
We now set out to bound the terms in the numerator of \cref{eq:limit} in terms of $\|e\|^2$:
\begin{itemize}
\item
For the gradient terms we get (by Cauchy-Schwartz):
\[
|(\nabla_x f(\xi,\overline{y}(x + e)) - \nabla_x f(x,\overline{y}(x)))\cdot e| \leq \|\nabla_x f(\xi,\overline{y}(x + e)) - \nabla_x f(x,\overline{y}(x))\|\cdot \|e\|.
\]
Furthermore, by Lipschitz continuinty of $\nabla_x f$ we have
\[
\|\nabla_x f(\xi,\overline{y}(x + e)) - \nabla_x f(x,\overline{y}(x))\| \leq L_{xx}\|\xi - x\| + L_{xy}\|\overline{y}(x + e) - \overline{y}(x)\|.
\]
Using that $\xi = x + t e$ and Lipschitz continuity of $\overline{y}$ (\cref{lma:ybar}) we find
\[
\|\nabla_x f(\xi,\overline{y}(x + e)) - \nabla_x f(x,\overline{y}(x))\| \leq \left(L_{xx} + L_{xy}L_{yx}/\mu)\|e\|\right).
\]
This results in
\[
|(\nabla_x f(\xi,\overline{y}(x + e)) - \nabla_x f(x,\overline{y}(x)))\cdot e| \leq (L_{xx} + L_{xy}L_{yx}/\mu)\|e\|^2.
\]
\item
Next, we need to bound $|F(x,\overline{y}(x+e)) - F(x,\overline{y}(x))|$ in terms of $\|e\|^2$. Since $\overline{y}(x)$ is the optimal solution at $x$ we have $|F(x,\overline{y}(x+e)) - F(x,\overline{y}(x))| = F(x,\overline{y}(x+e)) - F(x,\overline{y}(x))$.
Using the convexity of $F(x,\cdot)$ we have
\[
F(x,\overline{y}(x+e)) - F(x,\overline{y}(x)) \leq \langle u, \overline{y}(x+e)) - \overline{y}(x)\rangle,
\]
with $u \in \partial_y F(x,\overline{y}(x+e))$. This can be expressed as $u = \nabla_y f(x,\overline{y}(x+e)) + w$ with $w\in \partial r_2(\overline{y}(x+e))$. Choosing $w$ appropriately as $w = - \nabla_y f(x+e, \overline{y}(x+e))$ we get
\begin{multline*}
F(x,\overline{y}(x+e)) - F(x,\overline{y}(x)) \\
\leq \|\nabla_y f(x,\overline{y}(x+e)) - \nabla_y f(x + e,\overline{y}(x+e))\| \|\overline{y}(x+e)) - \overline{y}(x)\|.
\end{multline*}
Using Lipschitz continuity of $\nabla_y f$ and $\overline{y}$ we get
\[
F(x,\overline{y}(x+e)) - F(x,\overline{y}(x)) \leq (L_{yx}L_{xy}/\mu) \|e\|^2.
\]
\end{itemize}
We have thus upperbounded the term in \cref{eq:limit} in terms of $\|e\|$. As this upper bound tends to zero as $\|e\|\rightarrow 0$ and the fraction is always nonnegative we conclude that the limit tends to zero.
\end{proof}

\section{Proof of \cref{col:lipgrad}}
\begin{proof}
For ease of notation we fix $x,x'$ and denote the corresponding (unique) optimal solutions by $\overline{y}$ and $\overline{y}'$ respectively.
Using the definition $\nabla\overline{f}(x) = \nabla_x f(x,\overline{y})$ we have
\[
\|\nabla_x f(x,\overline{y}) - \nabla_x f(x',\overline{y}')\| \leq \|\nabla_x f(x,\overline{y}) - \nabla_x f(x',\overline{y})\| + \|\nabla_x f(x',\overline{y}) - \nabla_x f(x',\overline{y}')\|.
\]
Using Lipschitz continuity of $\nabla_x f$ and $\overline{y}$ we immediately find
\[
\|\nabla_x f(x,\overline{y}) - \nabla_x f(x',\overline{y}')\| \leq (L_{xx} + L_{xy}L_{yx}/\mu) \|x - x'\|
\]
\end{proof}

% \section{Gradient descent with errors}
% Consider solving
% \[
% \min_{x} f(x),
% \]
% with an inexact gradient method
% \[
% x_{k+1} = x_k - \alpha \nabla f(x_k) - \alpha e_k.
% \]
% Assume that $f$ is $L-$Lipschitz smooth, so we have
% \[
% f(y) \leq f(x) + \nabla f(x)^T(y - x) + \frac{L}{2}\|y - x\|^2,
% \]
% for any $x,y$. Substituting $x_{k+1}$ and $x_k$ we find
% \[
% f_{k+1} - f_k \leq -\alpha \nabla f_k^T(\nabla f_k + e_k) + \frac{\alpha^2L}{2}\|\nabla f_k + e_k\|^2.
% \]
% For $\alpha = 1/L$ we get a simplified expression
% \[
% f_{k+1} - f_k \leq \frac{1}{2L}\|e_k\|^2 - \frac{1}{2L}\|\nabla f_k\|^2,
% \]
% meaning that we can expect decrease of the objective as long as $\|e_k\|^2 < \|\nabla f_k\|^2$.
% Rearranging terms we get
% \[
% \|\nabla f_k\|^2 \leq 2L(f_k - f_{k+1}) + \|e_k\|^2.
% \]
% Taking a sum
% \[
% \sum_{k=0}^{n-1} \|\nabla f_k\|^2 \leq 2L(f_0 - f_{n-1}) + \sum_{k=0}^{n-1}\|e_k\|^2.
% \]
% Since $f_{n-1} \geq f_*$ we get
% \[
% \sum_{k=0}^{n-1} \|\nabla f_k\|^2 \leq 2L(f_0 - f_*) + \sum_{k=0}^{n-1}\|e_k\|^2.
% \]
% Thus if the errors are square-summable, so will the gradients implying that we have $\|\nabla f_k\|  = \mathcal{O}(1/k)$. Moreover,
% since $n \left(\min_{k\in {0,\ldots,n-1}} \|\nabla f_k\|^2\right) \leq \sum_{k=0}^{n-1} \|\nabla f_k\|^2$ we find
% \[
% \min_{k\in \{0,\ldots,n-1\}} \|\nabla f_k\|^2 \leq \frac{2L(f_0 - f_*)}{n} + \frac{1}{n}\sum_{k=0}^{n-1}\|e_k\|^2.
% \]
% Now we only need $\lim_{n\rightarrow\infty} \frac{1}{n}\sum_{k=0}^{n-1}\|e_k\|^2 < \infty$, implying a slower decay?

\section{Proximal gradient with errors}
\label{sec:inexactProofs}
Consider solving
\[
\min_x F(x),
\]
with $F(x) = f(x) + g(x)$, where $f$ is proper, closed and $L-$Lipschitz smooth and $g$ proper closed and convex. We consider an inexact proximal gradient method of the form
\begin{equation}
x_{k+1} = \text{prox}_{\alpha g} \left(x_k - \alpha(\nabla f(x_k) + e_k)\right).
\label{eq:basiciterappendix}
\end{equation}
For the following we closely follow \cite{Beck2017}.

\begin{lemma}[Sufficient decrease of inexact proximal gradient]\label{lemma:decrease}
% Introduce the prox-gradient mapping $T$ and the gradient mapping, $G(x)$:
% \[
% T(x) = \ts{\prox_{\alpha g}}(x - \alpha \nabla f(x))
% \]
% \[
% G(x) = \alpha^{-1}(x - T(x)).
% \]
The iteration \cref{eq:basiciterappendix} produces iterates that obey
\[
\begin{aligned}
F(x_k) - F(x_{k+1}) &\geq \left(\alpha^{-1} - L/2\right)\|x_{k+1} - x_{k}\|^2 - \langle e_k, x_{k+1} - x_{k} \rangle
%&\geq  \left(\alpha^{-1} - L/2\right)\|x_{k+1} - x_{k}\|^2 - \|e_k\|^2 -\alpha \|e_k\| \|G(x_k)\|  % \langle e_k, x_{k+1} - x_{k}
\end{aligned}
\]
\end{lemma}
\begin{proof}
By the smoothness assumption we have
\[
f(x_{k+1}) \leq f(x_{k}) + \langle \nabla f(x_k), x_{k+1} - x_{k}\rangle + \frac{L}{2}\|x_{k+1} - x_{k}\|^2.
\]
By \cite[Thm. 6.39]{Beck2017} we have
\[
\langle x_k - \alpha \nabla f(x_k) - \alpha e_k - x_{k+1}, x_k - x_{k+1} \rangle \leq \alpha g(x_k) - \alpha g(x_{k+1}),
\]
which yields
\[
\langle \nabla f(x_k), x_{k+1} - x_k\rangle \leq - \alpha^{-1}\|x_{k+1} - x_{k}\|^2 - \langle e_k, x_{k+1} - x_{k} \rangle + g(x_k) - g(x_{k-1}).
\]
Combining gives the result.
%To get the second inequality, we observe that
%\[
%\begin{aligned}
%\langle e_k, x_{k+1} - x_{k} \rangle  & = \langle e_k, \ts{\prox_{\alpha g}}(x_k - \alpha \nabla f(x_k)-\alpha e_k) - x_k\rangle  \\
%& = \langle e_k, \ts{\prox_{\alpha g}}(x_k - \alpha \nabla f(x_k)-\alpha e_k)  - T(x_k)\rangle + \langle e_k, T(x_k) - x_k\rangle
%\end{aligned}
%\]
%The prox map is non-expansive, and so we have the simple bound
%\[
%|\langle e_k, x_{k+1} - x_{k} \rangle | \leq \|e_k\|^2 + \alpha\|e_k\| \|G(x_k)\|.
%\]
\end{proof}

We let $\alpha < 2/L$ and state a simple corollary that ensures that we can get descent if $\|e_k\|$ is small enough.

\begin{corollary}[Existence of small enough errors]\label{cor:decrease}
At any iterate, we can always take $\|e_k\|$ small enough to ensure descent, unless $x_k$ is a stationary point.
\end{corollary}
\begin{proof}
Introduce the prox-gradient mapping $T$:
\[
T(x) = \ts{\prox_{\alpha g}}(x - \alpha \nabla f(x))
\]
If $x_k$ is non-stationary, we know that $\gamma := \|T(x_k) - x_k\|$ is bounded away from $0$.
On the  other hand, the function
\[
h(e_k) :=\|\prox_{\alpha g} (x_k - \alpha \nabla f(x_k)-\alpha e_k) - x_k\|
\]
is $\alpha$-Lipschitz continuous, since the norm and the prox map are both 1-Lipschitz continuous,
and we have $h(0)  = \|T(x_k)-x_k\|$ and $h(e_k) = \|x_{k+1} - x_k\|$.
Therefore, if we take $\|e_k\| \leq \frac{\gamma}{2\alpha}$,
we have
\[
\| h(0) - h(e_k)\| \leq   \alpha \|e_k\|  \leq \frac{\gamma}{2}
\]
and therefore
\[
| \| x_{k+1} - x_k\| - \gamma| < \frac{\gamma}{2} \Rightarrow  \| x_{k+1} - x_k\|  > \frac{\gamma}{2}.
\]
So if we take $\|e_k\| \leq \min\left(\frac{\gamma}{2\alpha},  \frac{\gamma}{(2 (\alpha^{-1} - L/2)}\right)$,
we are guaranteed that $\|e_k\| < (\alpha^{-1} - L/2)\|x_{k+1} - x_k\|$, ensuring descent by \cref{lemma:decrease}.
\end{proof}

%We let $\alpha < 2/L$. The result implies that we need $\|e_k\| < (\alpha^{-1} - L/2)\|x_{k+1} - x_k\|$ to get descent at each iteration.
%This is guaranteed to hold for all small $e_k$ by  the non-expansiveness (1-Lipschitz continuity) of the prox map.
%Specifically, as $\|e_k\| \downarrow 0$, we have $(x_k - \alpha \nabla f(x_k)-\alpha e_k)  \rightarrow (x_k - \alpha \nabla f(x_k))$,
%and so $x_{k+1} \rightarrow T(x_k)$, and $\|x_{k+1} -x_k\| \rightarrow \|T(x_k) - x_k\|$
%which is positive unless $x_k$ is a stationary point. Supposing
%\[
% \|T(x_k) - x_k\| = \eta > 0,
%\]
%by Lipschitz continuity of the

%It seems that we cannot guarantee this at all in general; we can think of an adverserial error for which $x_{k+1} = x_k$. However, making such an error smaller should eventually lead to the desired bound.

%\textcolor{red}{Need to show that for any $x_k$ that is not a fixed point, there exists an $e$ sufficiently small (but non-zero) such that $\|e_k\| \leq (\alpha^{-1} - L/2)\|x_{k+1} - x_k\|$}
%\subsection*{Proof of~\cref{eq:inexactproxgrad}}

\begin{theorem}[Convergence of inexact proximal gradient - general case]
\label{thm:inexactproxgradappendix}
The iteration \cref{eq:basiciterappendix} with stepsize $\alpha = 1/L$ and errors obeying $\|e_k\| \leq C \|x_{k+1} - x_k\|$ with $C < L/2$ produces iterates for which
\[
\min_{k\in\{0, 1, \ldots, n-1\}} \|x_{k+1} - x_k\| \leq A\sqrt{\frac{F(x_0) - F(x_*)}{n}},
\]
with $A=\sqrt{\frac{2}{L-2C}}$.
\end{theorem}
From \cref{lemma:decrease} we have
\[
\|x_{k+1} - x_k\|^2 \leq \frac{2}{L}\left((F(x_{k}) - F(x_{k+1})) + \cos \theta_k \|e_k\|\|x_{k+1} - x_k\|\right)
\]
Summing over $k$ we have
\[
\sum_{k=0}^{n-1} \|x_{k+1} - x_k\|^2 \leq \frac{2}{L}(F(x_0) - F(x_*)) + \frac{2}{L}\sum_{k=0}^{n-1} \|e_k\|\|x_{k+1} - x_k\|
\]
Assuming that $\|e_k\| \leq C \|x_{k+1} - x_k\|$ we get
\[
\sum_{k=0}^{n-1} \|x_{k+1} - x_k\|^2 \leq  \frac{2}{L}(F(x_0) - F(x_*)) + \frac{2C}{L}\sum_{k=0}^{n-1} \|x_{k+1} - x_k\|^2,
\]
so
\[
\left(1 - \frac{2C}{L}\right) \sum_{k=0}^{n-1} \|x_{k+1} - x_k\|^2 \leq \frac{2}{L}(F(x_0) - F(x_*)).
\]
hence
\[
\sum_{k=0}^{n-1} \|x_{k+1} - x_k\|^2 \leq  \left(1 - \frac{2C}{L}\right)^{-1}\frac{2}{L} (F(x_0) - F(x_*))
\]
thus
\[
\min_{k\in\{0, 1, \ldots, n-1\}} \|x_{k+1} - x_k\| \leq \sqrt{\frac{2(F(x_0) - F(x_*))}{(L-2C)n}}
\]

\section{A stopping criterion for prox gradient descent}
When solving
\[
\min_x f(x) + g(x),
\]
where $f$ is $L-$Lipschitz-smooth and $f + g$ is $\mu-$strongly convex, with a proximal gradient method, we need a practical stopping criterion that can guarantee a certain distance to the optimal solution. The following bound is usefull:

\begin{lemma}[A stopping criterion for proximal gradient descent]
\label{lma:stopping}
Introduce the prox-gradient mapping $T$ and the gradient mapping, $G(x)$:
\[
T(x) = \ts{\prox_{\alpha g}}(x - \alpha \nabla f(x))
\]
\[
G(x) = \alpha^{-1}(x - T(x)),
\]
where $f$ is $L-$Lipschitz-smooth and $\mu-$strongly convex, $g$ is convex and $\alpha \in (0,2/L)$. Define the fixed point of $G$ by $\overline{x}$. Then,
\[
\|T(x) - \overline{x}\| \leq \frac{1 + \alpha L }{\mu}\|G(x)\|.
\]
\end{lemma}
\begin{proof}
By strong convexity of $F = f + g$ we have
\[
\mu \|T(x) - \overline{x}\|^2 \leq \langle d, T(x) - \overline{x} \rangle,
\]
with $d \in \partial F(T(x))$. Note that $\partial F(x) = \nabla f(x) + \partial g(x)$ and that
\[
\alpha^{-1}(x - T(x)) - \nabla f(x) \in \partial g(T(x)),
\]
by definition of the proximal operator. Picking this particular element in $\partial g(T(x))$ we obtain
\[
\mu \|T(x) - \overline{x}\|^2 \leq \langle G(x) ,T(x) - \overline{x} \rangle + \langle \nabla f(T(x)) - \nabla f(x) ,T(x) - \overline{x} \rangle
\]
By Lipschitz-smoothness and Cauchy-Schwartz we get
\[
\|T(x) - \overline{x}\| \leq \frac{1 + \alpha L }{\mu}\|G(x)\|.
\]
\end{proof}
\begin{remark}
When applying a standard proximal gradient method $x^+ = T(x)$ we can immediately use this bound to devise a stopping criterion $\|x^+ - x\|\leq \epsilon$ that guarantees $\|x^+ - \overline{x}\| \leq \epsilon (\alpha^{-1} + L) /\mu$.
\end{remark}


\begin{thebibliography}{10}

\bibitem{aravkin2016smart}
Aleksandr Aravkin and Damek Davis.
\newblock A smart stochastic algorithm for nonconvex optimization with
  applications to robust machine learning.
\newblock {\em arXiv preprint arXiv:1610.01101}, 2016.

\bibitem{AravkinVanLeeuwen2012}
Aleksandr~Y. Aravkin and Tristan van Leeuwen.
\newblock Estimating nuisance parameters in inverse problems.
\newblock {\em Inverse Problems}, 28(11):115016, nov 2012.

\bibitem{Austin2019}
Anthony~P. Austin, Zichao Wendy, Sven Leyffer, and Stefan~M. Wild.
\newblock {Simultaneous Sensing Error Recovery and Tomographic Inversion Using
  an Optimization-Based Approach}.
\newblock {\em SIAM Journal on Scientific Computing}, 41(3):B497--B521, jan
  2019.

\bibitem{Beck2017}
Amir Beck.
\newblock {\em {First-Order Methods in Optimization}}.
\newblock Society for Industrial and Applied Mathematics, 2017.

\bibitem{Bell2008ADo}
B.M. Bell and J.V. Burke.
\newblock Algorithmic differentiation of implicit functions and optimal values.
\newblock In Christian~H. Bischof, H.~Martin B{\"u}cker, Paul~D. Hovland, Uwe
  Naumann, and J.~Utke, editors, {\em Advances in Automatic Differentiation},
  pages 67--77. Springer, 2008.

\bibitem{Cornelio2012}
Anastasia Cornelio, E.~Loli Piccolomini, and J.~G. Nagy.
\newblock {Constrained variable projection method for blind deconvolution}.
\newblock In {\em Journal of Physics: Conference Series}, volume 386, pages
  6--11, 2012.

\bibitem{devolder2014first}
Olivier Devolder, Fran{\c{c}}ois Glineur, and Yurii Nesterov.
\newblock First-order methods of smooth convex optimization with inexact
  oracle.
\newblock {\em Mathematical Programming}, 146(1-2):37--75, 2014.

\bibitem{GolubPereyra}
{G.H.} Golub and {V.} Pereyra.
\newblock The differentiation of pseudo-inverses and nonlinear least squares
  which variables separate.
\newblock {\em SIAM J. Numer. Anal.}, 10(2):413--432, 1973.

\bibitem{GolubPereyra2003}
G.H. Golub and V.~Pereyra.
\newblock Separable nonlinear least squares: the variable projection method and
  its applications.
\newblock {\em Inverse Problems}, 19(2):R1, 2003.

\bibitem{Osborne2007}
M.R. Osborne.
\newblock Separable least squares, variable projection, and the
  {G}auss-{N}ewton algorithm.
\newblock {\em Electronic Transactions on Numerical Analysis}, 28(2):1--15,
  2007.

\bibitem{Pereyra2012}
V.~Pereyra and G.~Scherer, editors.
\newblock {\em {Exponential Data Fitting and its Applications}}.
\newblock {B}entham {S}cience and {S}cience {P}ublishers, mar 2012.

\bibitem{RTRW}
R.T. Rockafellar and R.J.B. Wets.
\newblock {\em Variational Analysis}, volume 317.
\newblock Springer, 1998.

\bibitem{rousseeuw1984least}
P.~J Rousseeuw.
\newblock {Least Median of Squares Regression}.
\newblock {\em Journal of the American statistical association},
  79(388):871--880, 1984.

\bibitem{ruhe1980algorithms}
A.~Ruhe and P.A. Wedin.
\newblock Algorithms for separable nonlinear least squares problems.
\newblock {\em SIAM Rev.}, 22(3):318--337, 1980.

\bibitem{schmidt2011convergence}
Mark Schmidt, Nicolas~L Roux, and Francis~R Bach.
\newblock Convergence rates of inexact proximal-gradient methods for convex
  optimization.
\newblock In {\em Advances in neural information processing systems}, pages
  1458--1466, 2011.

\bibitem{Shearer2013}
Paul Shearer and Anna~C Gilbert.
\newblock {A generalization of variable elimination for separable inverse
  problems beyond least squares}.
\newblock {\em Inverse Problems}, 29(4):045003, apr 2013.

\bibitem{Smith1992}
Ronald~L. Smith.
\newblock {Some interlacing properties of the Schur complement of a Hermitian
  matrix}.
\newblock {\em Linear Algebra and its Applications}, 177:137--144, 1992.

\bibitem{van2017automatic}
Tristan van Leeuwen, Simon Maretzke, and K~Joost Batenburg.
\newblock {Automatic alignment for three-dimensional tomographic
  reconstruction}.
\newblock {\em Inverse Problems}, 34(2):024004, feb 2018.

\bibitem{wang2015projection}
Weiran Wang and Canyi Lu.
\newblock Projection onto the capped simplex.
\newblock {\em arXiv preprint arXiv:1503.01002}, 2015.

\bibitem{yang2018general}
Eunho Yang, Aur{\'e}lie~C Lozano, Aleksandr Aravkin, et~al.
\newblock A general family of trimmed estimators for robust high-dimensional
  data analysis.
\newblock {\em Electronic Journal of Statistics}, 12(2):3519--3553, 2018.

\end{thebibliography}
\end{document}